%% file: metal.tex
\documentclass[12pt,unicode]{amsart}

\usepackage{cmap}
\usepackage[cp1251]{inputenc}
\usepackage[TS1,T2A]{fontenc}
\usepackage[russian]{babel}

\usepackage{geometry}
\usepackage{amsmath}
\usepackage{amsfonts}
\usepackage{amssymb}

\usepackage{pgf,tikz,circuitikz}
\usetikzlibrary{arrows}
\usepackage{graphicx}

\usepackage[pdftex,%
            colorlinks,linkcolor=red,citecolor=blue,urlcolor=blue,%
            pdfstartview=FitH,%
            pdfauthor={Sergey Dorichenko, Maxim Prasolov, Mikhail Skopenkov},%
            pdftitle={Dissections of a metal rectangle},%
            pdfkeywords={Dissection, Rectangle, Electrical network.}%
           ]{hyperref}

\newcommand{\mscomm}[1]{}

\theoremstyle{definition}
\newtheorem{pr}{Задача}

\begin{document}


\title{Квадрат из подобных прямоугольников}

\author{С. Дориченко и М. Скопенков}

\address{Школа 57 города Москвы}
\email{sdorichenko@gmail.com}

\address{Институт проблем передачи информации Российской Академии Наук} 
\email{skopenkov@rambler.ru}

\date{}

\maketitle

\newcounter{problem}
\setcounter{problem}0




\footnotetext[0]{Работа выполнена при поддержке гранта Президента РФ
МК-3965.2012.1, фонда <<Династия>>, фонда Саймонса}

\mscomm{
\bigskip
\textbf{План.}

1. Это статья - продолжение. Рекомендуется прочитать часть I, а если нет - то в некоторые утверждения придется поверить.

2. Постановки задачи

3. Мы решим эту задачу в частном случае, а в общем приведем только ответ.

4. Если отношение сторон прямоугольника рационально, то квадрат можно разрезать на прямоугольники, подобные ему

5. Приведем пример, коглда иррационально. Рассмотрим картинку с 5 прямоугольниками. Предположим, что они подобны и все найдем.

6. Задачи с тремя прямоугольниками

7. Попробуем обобщить пример, который мы разобрали. Одни укладываем так, другие сяк. 

8. Решаем это уравнение. Получаем описание всех отношений, которые могут получиться при таком способе разрезаний. 

9. Задача - разрежьте на прямоугольники с отношением $2+\sqrt{2}$.

10. мы показали, что при нашем способе разрезания ничего другого получить нельзя. Покажем, что другие способы разрезания ...

11. Напоминаие о том, что было в 1-й части

12. Алгебраичность 

13. Формулировка ключевого примера

14. Определение хорошего числа

15. Простейшие алгебраические свойства

16. Удивительная лемма

17. Доказательство ключевого примера

18. Ответ для хороших чисел

19. Доказательство этого ответа

20. Формулировка ответа на исходный вопрос

21. Какие-то эквивалентности в этой теореме

22. Про переменный ток, очень неформально

23. Задачки

24. Решения к упражненям и задачам

\bigskip}

Это статья --- продолжение статьи ``Разрезание металлического прямоугольника'' из  ``Кванта''  N3 за 2011 г. Для понимания этой статьи желательно прочитать предыдущую, но можно читать эту статью и независимо; правда, тогда придется поверить в некоторые утверждения.

В предыдущей статье мы ответили на вопрос, какие прямоугольники можно разрезать на квадраты (не обязательно одинакового размера). 
Вопрос, который нас интересует теперь, таков:
{\it когда из прямоугольников, подобных данному, можно сложить квадрат?} Более точно, каким может быть отношение сторон такого прямоугольника? Мы решим эту задачу в частном случае, а в общем приведем только ответ.





Разумеется, если стороны прямоугольника относятся как $m : n$, где $m$ и $n$ --- целые, то тогда можно; см. рис.~\ref{ris0d}.  На нашем рисунке все прямоугольники равны и расположены ``одинаково''. Ясно, что при таком способе укладки нам удастся сложить квадрат, только если отношение сторон рационально.

\begin{figure}[h]
\input{metal-fig0d.tex}
\caption{Разрезание квадрата на подобные прямоугольники с отношением сторон $m:n$}\label{ris0d}
\end{figure}

\medskip
\textbf{Может ли $R$ быть иррациональным}
\smallskip



А если класть некоторые прямоугольники поперек? Например,
посмотрите на рис.~\ref{ris4} --- могут ли все прямоугольники на нем оказаться подобными? 

Предположим, что это так: пусть все прямоугольники на этом рисунке имеют одно и то же отношение сторон $R$. 
Пусть сторона квадрата равна $1$. Тогда последовательно находим
$AB=1/3$, $AC=R/3$, $CD=1-R/3$, $DE=R-R^2/3$. С другой стороны, $DE=1/2$, значит, $R-R^2/3=1/2$. Решая квадратное уравнение,
находим $R=(3\pm\sqrt{3})/2$~--- иррациональное число.

Нетрудно убедиться, что при $R=(3\pm\sqrt{3})/2$ такое разрезание, как на рис.~\ref{ris4}, действительно возможно.
\mscomm{\textbf{Картинка действительно возможна при таких $R$}}

\smallskip

\begin{figure}[h]
\input{metal-fig7.tex}
\caption{Разрезание квадрата на $5$ подобных прямоугольников}\label{ris4}
\end{figure}

\bigskip
\bigskip

{
\small

\begin{pr} \label{p-9} Дизайнеру заказали раму для квадратного окна.
На проекте (рисунок~\ref{ris9})
показано, как должны примыкать стекла друг к другу, но их размеры искажены. Можно ли сделать все стекла в раме подобными прямоугольниками? 
\end{pr}

}

\begin{figure}[h]
\input{metal-fig7a.tex} 
    \caption{Проект оконной рамы.} 
    \label{ris9}
\end{figure}

Попробуем обобщить пример, который мы разобрали.

{\small
\begin{pr}\label{vdol-poperek}
Имеется два набора прямоугольников с отношением сторон $R$, в каждом из которых все прямоугольники равны. 
Прямоугольники первого набора сложили в один прямоугольник так, что получилось $m$ рядов по $n$ прямоугольников. Прямоугольники второго набора сложили в один прямоугольник так, что получилось $p$ рядов по $q$ прямоугольников.  После этого полученные два ``составных'' прямоугольника приложили друг к другу так, что прямоугольники первого набора оказались расположены ``вдоль'', а второго --- ``поперек''. Известно, что в результате получился квадрат. Найдите, какие значения может принимать число $R$ при данных $m,n,p,q$.
\end{pr}

\begin{pr} Разрежьте квадрат на несколько прямоугольников с отношениями сторон\\ \textbf{(A)} $2+\sqrt{2}$; \textbf{(B)} $2-\sqrt{2}$. 
\end{pr}

\begin{pr} \footnote{А. Шаповалов, Турнир городов} Можно ли разрезать квадрат на $3$ неравных подобных прямоугольника?
\end{pr}

}





\medskip

{\bf Сопротивление электрической цепи}

\smallskip

Нам снова поможет физическая интерпретация. В конце предыдущей статьи мы обещали подробно рассказать о сопротивлении электрической цепи. Выполним это обещание.


Вернемся к математической модели электрической цепи --- графу, ребрам которого сопоставлены положительные числа (\emph{сопротивления} ребер и \emph{напряжение} батарейки). Напомним, что \emph{силы токов} через ребра --- это сопоставленные ребрам действительные числа, которые определяются правилами Кирхгофа.
В дальнейшем будем считать, что напряжение батарейки равно $1$. Тогда назовем \emph{сопротивлением} электрической цепи величину, обратную силе тока через батарейку.

В предыдущей статье мы доказали, что для прямоугольника с горизонтальной стороной $1$, разрезанного на квадраты, токи в соответствующей электрической цепи равны длинам сторон квадратов. Однако наши рассуждения почти дословно переносятся на случай разрезания на прямоугольники.
Положим сопротивления резисторов равными отношениям сторон соответствующих прямоугольников. Тогда силы тока в резисторах будут равны длинам вертикальных сторон меньших прямоугольников, а сила тока через батарейку --- длине вертикальной стороны большого. Значит, \emph{сопротивление цепи равно отношению сторон большого прямоугольника}. (Мы уже объяснили это в первой части статьи на ``физическом уровне строгости'', а сейчас доказали математически строго.)




\begin{figure}[h]
\begin{tabular}[t]{c}
$R=R_1+R_2$\\[3pt]
\input{metal-fig5a.tex}
\end{tabular}
\begin{tabular}[t]{c}
$R=\frac{R_1R_2}{R_1+R_2}$\\[2pt]
\input{metal-fig5b.tex}\\
\end{tabular}
\caption{
Формулы для сопротивления цепей из последовательно (слева) и параллельно (справа) соединенных резисторов}\label{risposlpar}
\end{figure}

Как же найти сопротивление цепи? Известны формулы для сопротивления цепей из последовательно и параллельно соединенных резисторов; см.
рисунок~\ref{risposlpar}
(докажите эти формулы!). Оказывается, для произвольной электрической цепи тоже можно написать формулу, которая выражает сопротивление цепи через сопротивления отдельных резисторов:




\smallskip

\noindent{\bf Теорема о сопротивлении цепи.}
{\it 
Силы тока в электрической цепи и ее сопротивление можно выразить через сопротивления отдельных резисторов, используя только операции сложения, вычитания, умножения и деления.}

\smallskip

\noindent\textbf{Доказательство.} Мы докажем эту теорему и в дальнейшем будем использовать ее только для случая электрических цепей, построенных по разрезанию прямоугольника на прямоугольники\footnote{Также мы не доказываем и не используем, что сопротивление цепи выражается через сопротивления резисторов \emph{одинаковым образом} для любых значений сопротивлений резисторов.}.
Запишем систему уравнений на силы токов, построенную по правилам Кирхгофа. Ее коэффициенты --- это $\pm 1$, а также сопротивления отдельных резисторов.
Вспомним теорему о решении системы из предыдущей статьи: мы фактически доказали тогда, что если решение системы единственно, то оно выражается через коэффициенты с помощью четырех арифметических операций.
У нашей системы заведомо есть решение, потому что ей удовлетворяют длины вертикальных сторон прямоугольников разрезания. Других решений нет по теореме единственности. Значит, все силы тока выражается через $\pm 1$ и сопротивления резисторов 
с помощью арифметических операций. То же самое верно и для сопротивления цепи. Осталось заметить, что единицу можно представить в виде $1=R/R$, где $R$ --- сопротивление одного из резисторов, поэтому
теорема о сопротивлении цепи доказана.

\medskip

\medskip
\textbf{Может ли $R$ быть любым}
\smallskip


Может быть, тогда любое число $R$ годится? Оказывается, нет:


\smallskip

\noindent{\bf Теорема 1.} {\it Пусть квадрат разрезан на подобные прямоугольники, у которых отношение большей стороны к меньшей равно $R$. 
Тогда $R$~--- корень ненулевого многочлена с целыми коэффициентами.}

\smallskip


\noindent{\bf Доказательство}. Растянем нашу картинку в $R$ раз по горизонтали. Получим прямоугольник с отношением сторон $R$, разрезанный на квадраты и прямоугольники с отношением
сторон $R^2$. Построим соответствующую электрическую цепь.
Она состоит из резисторов сопротивлением $1$ и $R^2$, а сама имеет сопротивление $R$.

По теореме о сопротивлении цепи число $R$ можно выразить через числа $1$ и $R^2$, пользуясь только четырьмя арифметическими действиями. Это означает, что найдутся два ненулевых многочлена $p(x)$ и $q(x)$ с целыми коэффициентами, такие что
$R=\frac{p(R^2)}{q(R^2)}$. Значит, $q(R^2)\cdot R-p(R^2)=0$. Получаем, что число $R$~--- корень многочлена $q(x^2)\cdot x-p(x^2)$. Этот многочлен имеет целые коэффициенты. Он ненулевой, так как
многочлены $p(x^2)$ и $q(x^2)$~--- ненулевые, причём в многочлен $p(x^2)\cdot x$
переменная $x$ входит только в нечётной степени, а в многочлен $q(x^2)$~---
только в чётной. Теорема~1 доказана.

\smallskip

Но эта теорема не дает полного ответа на вопрос о том, какие отношения $R$ возможны.

\medskip
\textbf{Ключевой пример}
\smallskip

Существуют числа $R$, являющиеся корнями многочленов с целыми коэффициентами, для которых требуемое разрезание невозможно. Приведем ключевой пример такого числа.

\smallskip

\noindent{\bf Теорема 2.} {\it Квадрат нельзя разрезать на подобные прямоугольники, у которых отношение большей стороны к меньшей равно $1+\sqrt2$.}

\smallskip

Эту теоремы так сразу не докажешь, потребуется некоторая подготовка. Рассмотрим все числа, которые можно представить в виде $x=a+b\sqrt{2}$, где $a$ и $b$ --- рациональные числа. Числа такого вида назовем \emph{хорошими}. В дальнейшем, в выражении $a+b\sqrt{2}$ мы всегда подразумеваем, что $a$ и $b$ рациональны. Для числа $x=a+b\sqrt{2}$ назовем \emph{сопряженным} к нему число $\bar x=a-b\sqrt{2}$.

\smallskip

\noindent\textbf{Хорошие свойства хороших чисел.}
\begin{enumerate}
\item Если $x$ и $y$ --- хорошие, то и $x+y$, $x-y$, $x\cdot y$, $x/y$ --- хорошие.
\item Хорошее число $x$ единственным образом представляется в виде $x=a+b\sqrt{2}$.
\item Если $x$ и $y$ --- хорошие, то $\bar x+\bar y=\overline{x+y}$ и $\bar x\cdot\bar y=\overline{x\cdot y}$.
\end{enumerate}
\smallskip

\mscomm{Перед тем, как читать дальше,} Рекомендуем читателю самому попробовать доказать эти свойства.

\smallskip

Докажем удивительную лемму, в которой физика (где все числа заданы с ограниченной точностью) переплетается с алгеброй (где встречаются точные выражения вида $a+b\sqrt{2}$).

\smallskip

\noindent{\bf Удивительная лемма.} {\it Пусть в электрической цепи сопротивления всех резисторов --- хорошие числа, причем сопряженные к ним числа положительны. Тогда если заменить сопротивления резисторов на сопряженные к ним числа, то сопротивление цепи также заменится на сопряженное.}

\smallskip

\noindent\textbf{Доказательство.} 
По теореме о сопротивлении цепи все силы тока в исходной цепи
хорошие.
Заменим сопротивления всех резисторов и все силы тока на сопряженные. 
Из свойств хороших чисел следует, что правила Кирхгофа останутся справедливыми. Так как все новые сопротивления положительны, то выполнена теорема единственности.  Значит, выбранные силы тока и будут силами тока в полученной цепи.
Так как ток через батарейку заменился на сопряженный, то
и сопротивление цепи заменилось на сопряженное. Удивительная лемма доказана.

\smallskip

\noindent\textbf{Доказательство теоремы 2.} 
Растянем картинку в $1+\sqrt2$ раз по горизонтали. Получим разрезание прямоугольника с отношением сторон $1+\sqrt2$. Рассмотрим соответствующую электрическую цепь. Она состоит из 
резисторов сопротивлением $1$ и $(1+\sqrt2)^2$, а сопротивление цепи равно $1+\sqrt2$.
Это хорошие числа.

Заменим в нашей цепи сопротивления всех резисторов на сопряженные. Иными словами, заменим все резисторы сопротивлением $(1+\sqrt2)^2$ на резисторы сопротивлением $(1-\sqrt2)^2$. Сопротивления останутся положительными числами.
По удивительной лемме сопротивление цепи тоже заменится на сопряженное. Значит, оно станет равным $1-\sqrt{2}$. Но это число отрицательно!
Получаем противоречие (по задаче~5 из предыдущей статьи сила тока через батарейку, а следовательно, и сопротивление цепи, всегда положительны). Теорема 2 доказана.

\mscomm{
\medskip
\textbf{Основная теорема}
\smallskip
}

\smallskip


Мы готовы описать все хорошие числа $x$, для которых из прямоугольников с отношением сторон $x$ можно составить квадрат.

\smallskip

\noindent{\bf Теорема 3.} \emph{Пусть $x$ --- хорошее число. Тогда из прямоугольников с отношением сторон $x$ можно составить квадрат, если и только если оба числа $x$ и $\bar x$ положительны.}

\smallskip

\mscomm{Написать подробное доказательство!!!}
Доказательство этой теоремы мы оставляем в качестве задачи. Скажем лишь, что невозможность разрезания в случае $\bar x<0$ доказывается почти дословно так же, как теорема 2. Возможности разрезания при $x, \bar x>0$ доказывается с помощью конструкции из задачи~\ref{vdol-poperek}. 

\medskip
\textbf{Ответ на вопрос}
\smallskip

Теперь мы в шаге от ответа на вопрос, поставленный в начале статьи. 
В последней теореме все определялось знаком сопряженного числа. 
А в общей ситуации важен знак действительной части ``сопряженных'' чисел
\mscomm{Сноска: О комплексных числах можно прочитать в --- Дать ссылку!}:

\smallskip

\noindent{\bf Теорема Ласковича--Ринна--Секереша--Фрайлинга (1994)\footnote{Freiling C., Rinne D., Tiling a square with similar
rectangles, Math. Res. Lett. 1 (1994), 547--558; Laszkovich M., Szekeres G., Tiling of the square
with similar rectangles, Discr. Comp. Geometry 13 (1995),
569--572.}.} {\it Для числа $R>1$ следующие три условия эквивалентны:

\noindent 1) Квадрат можно разрезать на прямоугольники, у которых отношение большей стороны к меньшей равно $R$.

\noindent 2) Число $R$ является корнем некоторого многочлена с целыми коэффициентами, у которого все комплексные корни имеют положительную действительную часть.

\noindent 3) Для некоторых положительных рациональных чисел $c_i$ выполнено равенство}
$$c_1 R+\cfrac1{c_2 R +\cfrac1{c_3 R+\dots+\cfrac1 {c_n R}}}=1.$$

\smallskip

\mscomm{
\begin{figure}[ht]
\caption{Как разрезать квадрат на подобные прямоугольники.}
\end{figure}
}

\noindent{\bf Доказательство утверждения $3)  \Rightarrow 1)$.}
Пусть выполнено равенство из условия 3). Покажем, как разрезать квадрат на прямоугольники с отношением сторон $R$ и $1/R$. Возьмём некоторый квадрат. Отрежем от него прямоугольник с отношением сторон $c_1R$, проведя вертикальный разрез. От оставшейся части отрежем прямоугольник с отношением сторон $\frac{1}{c_2R}$, проведя горизонтальный разрез. Будем продолжать этот процесс, чередуя вертикальные и горизонтальные разрезы. В силу равенства в условии 3) на $n$-м шаге мы получим прямоугольник с отношением сторон $c_nR$ (или $\frac{1}{c_nR}$, в зависимости от чётности числа $n$). Тем самым квадрат оказался разбит на прямоугольники с отношениями сторон $c_1R$, $\frac{1}{c_2R}$, $c_3R$, $\dots$, $c_nR$ (или $\frac{1}{c_nR}$). Остаётся каждый из полученных прямоугольников разрезать на прямоугольники с отношением сторон $R$ или $1/R$, и нужное разбиение построено.

Доказательство остальных частей этой теоремы выходит за рамки настоящей статьи. Отметим только, что их можно доказать в целом аналогично теореме~2 выше, используя цепи \emph{переменного} тока\footnote{M.~Prasolov, M.~Skopenkov, Tiling by rectangles and alternating current, J. Combin. Theory Ser. A \textbf{118:3} (2011), 920--937, \url{http://arxiv.org/abs/1002.1356}.}, см. рис.~\ref{fig-perem}.

\begin{figure}[h]
\makebox[1.5cm][r]{
  \begin{tabular}{l}
  \input{metal-fig7d.tex}
  \end{tabular}
  } 
  \makebox[3.5cm][c]{
  \begin{tabular}{c}
  \input{circuit-figure.tex}
  \end{tabular}
  }
\caption{Физическая интерпретация разрезания на подобные прямоугольники, использующая цепи переменного тока.}
\label{fig-perem}
\end{figure}

В заключение приведем еще несколько задач.

\smallskip

\small

\begin{pr}\footnote{Вы можете найти решение в статье \mscomm{Дать ссылку на Квант! или}  S.~Wagon, Fourteen proofs of a result about tiling a rectangle, Amer.~Math.~Monthly 94:7 (1987).}
Прямоугольник разрезан на несколько прямоугольников, у каждого из которых длина хотя бы одной стороны --- целое число. Доказать, что и у исходного прямоугольника длина хотя бы одной стороны --- целое число.
\end{pr}

\begin{pr}
\footnote{А.~Шаповалов, Задача~7, XXVI Турнир Городов, \url{http://turgor.ru/26/index.php}}
Пусть $A$ и $B$ --- два прямоугольника. Из прямоугольников, равных $A$, сложили прямоугольник, подобный $B$.  Докажите, что из прямоугольников, равных $B$, можно сложить прямоугольник, подобный $A$.
\end{pr}

\begin{pr}\footnote{А.~Колотов, Об одном разбиении прямоугольника, Квант \No~1 (1973), \url{http://kvant.mccme.ru/1973/01/ob_odnom_razbienii_pryamougoln.htm}.}
Найдите необходимые и достаточные условия, которым должны удовлетворять числа $a$, $b$, $c$, $d$, чтобы прямоугольник $a\times b$ можно было разрезать на несколько прямоугольников $c\times d$.
\end{pr}


\normalsize

\mscomm{Написать решения задач!}

\end{document}

%% file: metal-fig0d.tex
\definecolor{zzttqq}{rgb}{0.6,0.2,0}
\begin{tikzpicture}[line cap=round,line join=round,>=triangle 45,x=0.3cm,y=0.3cm]
\clip(-3.34,-1.36) rectangle (4.1,5.96);
\fill[color=zzttqq,fill=zzttqq,fill opacity=0.1] (-3,5) -- (3,5) -- (3,-1) -- (-3,-1) -- cycle;
\draw [color=zzttqq] (-3,5)-- (3,5);
\draw [color=zzttqq] (3,5)-- (3,-1);
\draw [color=zzttqq] (3,-1)-- (-3,-1);
\draw [color=zzttqq] (-3,-1)-- (-3,5);
\draw [color=zzttqq] (-1.8,5)-- (-1.8,-1);
\draw [color=zzttqq] (-0.6,5)-- (-0.6,-1);
\draw [color=zzttqq] (0.6,5)-- (0.6,-1);
\draw [color=zzttqq] (1.8,5)-- (1.8,-1);
\draw [color=zzttqq] (-3,3)-- (3,3);
\draw [color=zzttqq] (-3,1)-- (3,1);
\begin{scriptsize}
\draw[color=zzttqq] (0.04,5.58) node {$m$};
\draw[color=zzttqq] (3.54,2.12) node {$n$};
\end{scriptsize}
\end{tikzpicture}

%% file: metal-fig7.tex
\definecolor{zzttqq}{rgb}{0.6,0.2,0}
\definecolor{ffqqqq}{rgb}{1,0,0}
\begin{tikzpicture}[line cap=round,line join=round,>=triangle 45,x=0.5cm,y=0.5cm]
\clip(-3.42,-2.1) rectangle (4.32,5.92);
\fill[color=zzttqq,fill=zzttqq,fill opacity=0.1] (-2,5) -- (4,5) -- (4,-1) -- (-2,-1) -- cycle;
\draw [color=zzttqq] (-2,5)-- (4,5);
\draw [color=zzttqq] (4,5)-- (4,-1);
\draw [color=zzttqq] (4,-1)-- (-2,-1);
\draw [color=zzttqq] (-2,-1)-- (-2,5);
\draw [color=zzttqq] (-2,3.58)-- (4,3.56);
\draw [color=zzttqq] (-0.02,3.57)-- (0,-1);
\draw [color=zzttqq] (2,3.57)-- (2,-1);
\draw [color=zzttqq] (1,5)-- (1,3.57);
\fill [color=ffqqqq] (-2,5) circle (1.5pt);
\draw[color=ffqqqq] (-2.78,5.46) node {$D$};
\fill [color=ffqqqq] (-2,-1) circle (1.5pt);
\draw[color=ffqqqq] (-2.72,-1.38) node {$A$};
\fill [color=ffqqqq] (-2,3.58) circle (1.5pt);
\draw[color=ffqqqq] (-2.72,3.74) node {$C$};
\fill [color=ffqqqq] (0,-1) circle (1.5pt);
\draw[color=ffqqqq] (0.56,-1.36) node {$B$};
\fill [color=ffqqqq] (1,5) circle (1.5pt);
\draw[color=ffqqqq] (1.12,5.48) node {$E$};
\end{tikzpicture}

%% file: metal-fig7a.tex
\definecolor{zzttqq}{rgb}{0.6,0.2,0}
\begin{tikzpicture}[line cap=round,line join=round,>=triangle 45,x=0.3cm,y=0.3cm]
\clip(-2.2,-1.16) rectangle (4.24,5.16);
\fill[color=zzttqq,fill=zzttqq,fill opacity=0.1] (-2,5) -- (4,5) -- (4,-1) -- (-2,-1) -- cycle;
\draw [color=zzttqq] (-2,5)-- (4,5);
\draw [color=zzttqq] (4,5)-- (4,-1);
\draw [color=zzttqq] (4,-1)-- (-2,-1);
\draw [color=zzttqq] (-2,-1)-- (-2,5);
\draw [color=zzttqq] (-2,3.58)-- (4,3.56);
\draw [color=zzttqq] (1.02,3.57)-- (1,-1);
\end{tikzpicture}

%% file: metal-fig5a.tex
\definecolor{zzttqq}{rgb}{0.6,0.2,0}
\begin{tikzpicture}[line cap=round,line join=round,>=triangle 45,x=0.5cm,y=0.5cm]
\clip(-3.14,-2.22) rectangle (2.22,5.14);
\fill[color=zzttqq,fill=zzttqq,fill opacity=0.1] (-3,5) -- (2,5) -- (2,1) -- (-3,1) -- cycle;
\fill[color=zzttqq,fill=zzttqq,fill opacity=0.1] (-2.62,-2) -- (-2.62,-1.46) -- (-0.74,-1.46) -- (-0.74,-2) -- cycle;
\fill[color=zzttqq,fill=zzttqq,fill opacity=0.1] (1.62,-2) -- (1.62,-1.46) -- (-0.26,-1.46) -- (-0.26,-2) -- cycle;
\draw [color=zzttqq] (-3,5)-- (2,5);
\draw [color=zzttqq] (2,5)-- (2,1);
\draw [color=zzttqq] (2,1)-- (-3,1);
\draw [color=zzttqq] (-3,1)-- (-3,5);
\draw [color=zzttqq] (-1,5)-- (-1,1);
\draw [color=zzttqq] (-2.62,-2)-- (-2.62,-1.46);
\draw [color=zzttqq] (-2.62,-1.46)-- (-0.74,-1.46);
\draw [color=zzttqq] (-0.74,-1.46)-- (-0.74,-2);
\draw [color=zzttqq] (-0.74,-2)-- (-2.62,-2);
\draw [color=zzttqq] (1.62,-2)-- (1.62,-1.46);
\draw [color=zzttqq] (1.62,-1.46)-- (-0.26,-1.46);
\draw [color=zzttqq] (-0.26,-1.46)-- (-0.26,-2);
\draw [color=zzttqq] (-0.26,-2)-- (1.62,-2);
\draw (-3,-1.74)-- (-2.62,-1.74);
\draw [color=zzttqq] (-0.74,-1.72)-- (-0.26,-1.72);
\draw [color=zzttqq] (1.62,-1.72)-- (2,-1.72);
\draw [color=zzttqq] (-3,-1.74)-- (-3,-0.5);
\draw [color=zzttqq] (2,-1.72)-- (2,-0.5);
\draw [color=zzttqq] (-3,-0.5)-- (-0.54,-0.5);
\draw [color=zzttqq] (2,-0.5)-- (-0.26,-0.5);
\draw [color=zzttqq] (-0.54,-0.02)-- (-0.54,-1);
\draw [line width=2.4pt,color=zzttqq] (-0.26,-0.24)-- (-0.26,-0.82);
\draw[color=zzttqq] (0.6,3.2) node {$R_2$};
\draw[color=zzttqq] (-1.86,3.18) node {$R_1$};
\draw[color=zzttqq] (-1.7,-1.06) node {$R_1$};
\draw[color=zzttqq] (0.84,-1.02) node {$R_2$};
\end{tikzpicture}

%% file: metal-fig5b.tex
\definecolor{zzttqq}{rgb}{0.6,0.2,0}
\begin{tikzpicture}[line cap=round,line join=round,>=triangle 45,x=0.5cm,y=0.5cm]
\clip(-3.14,-4.02) rectangle (2.22,5.14);
\fill[color=zzttqq,fill=zzttqq,fill opacity=0.1] (-3,5) -- (2,5) -- (2,1) -- (-3,1) -- cycle;
\fill[color=zzttqq,fill=zzttqq,fill opacity=0.1] (-1.44,-2) -- (-1.44,-1.52) -- (0.52,-1.52) -- (0.52,-2) -- cycle;
\fill[color=zzttqq,fill=zzttqq,fill opacity=0.1] (-1.44,-3) -- (-1.44,-2.52) -- (0.52,-2.52) -- (0.52,-3) -- cycle;
\draw [color=zzttqq] (-3,5)-- (2,5);
\draw [color=zzttqq] (2,5)-- (2,1);
\draw [color=zzttqq] (2,1)-- (-3,1);
\draw [color=zzttqq] (-3,1)-- (-3,5);
\draw [color=zzttqq] (-1.44,-2)-- (-1.44,-1.52);
\draw [color=zzttqq] (-1.44,-1.52)-- (0.52,-1.52);
\draw [color=zzttqq] (0.52,-1.52)-- (0.52,-2);
\draw [color=zzttqq] (0.52,-2)-- (-1.44,-2);
\draw [color=zzttqq] (-3,-1.74)-- (-1.44,-1.74);
\draw [color=zzttqq] (-3,-1.74)-- (-3,-0.5);
\draw [color=zzttqq] (2,-1.72)-- (2,-0.5);
\draw [color=zzttqq] (-3,-0.5)-- (-0.54,-0.5);
\draw [color=zzttqq] (2,-0.5)-- (-0.26,-0.5);
\draw [color=zzttqq] (-0.54,-0.02)-- (-0.54,-1);
\draw [line width=2.4pt,color=zzttqq] (-0.26,-0.24)-- (-0.26,-0.82);
\draw [color=zzttqq] (-3,3)-- (2,3);
\draw [color=zzttqq] (2,-1.72)-- (0.52,-1.72);
\draw [color=zzttqq] (0.52,-3)-- (-1.44,-3);
\draw [color=zzttqq] (0.52,-2.52)-- (0.52,-3);
\draw [color=zzttqq] (-1.44,-2.52)-- (0.52,-2.52);
\draw [color=zzttqq] (-1.44,-3)-- (-1.44,-2.52);
\draw [color=zzttqq] (-3,-2.74)-- (-1.44,-2.74);
\draw [color=zzttqq] (2,-2.72)-- (0.52,-2.72);
\draw [color=zzttqq] (-1.44,-3)-- (-1.44,-2.52);
\draw [color=zzttqq] (-1.44,-2.52)-- (0.52,-2.52);
\draw [color=zzttqq] (0.52,-2.52)-- (0.52,-3);
\draw [color=zzttqq] (0.52,-3)-- (-1.44,-3);
\draw [color=zzttqq] (-3,-1.74)-- (-3,-2.74);
\draw [color=zzttqq] (2,-1.72)-- (2,-2.72);
\draw[color=zzttqq] (-0.36,2.14) node {$R_2$};
\draw[color=zzttqq] (-1.58,-1.02) node {$R_1$};
\draw[color=zzttqq] (-0.38,4.04) node {$R_1$};
\draw[color=zzttqq] (1.1,-3.52) node {$R_2$};
\end{tikzpicture}

%% file: metal-fig7d.tex
\definecolor{zzttqq}{rgb}{0.6,0.2,0}
\begin{tikzpicture}[line cap=round,line join=round,>=triangle 45,x=0.3cm,y=0.3cm]
\clip(-3.42,-2.1) rectangle (4.32,5.92);
\fill[color=zzttqq,fill=zzttqq,fill opacity=0.1] (-2,5) -- (4,5) -- (4,-1) -- (-2,-1) -- cycle;
\draw [color=zzttqq] (-2,5)-- (4,5);
\draw [color=zzttqq] (4,5)-- (4,-1);
\draw [color=zzttqq] (4,-1)-- (-2,-1);
\draw [color=zzttqq] (-2,-1)-- (-2,5);
\draw [color=zzttqq] (-2,3.58)-- (4,3.56);
\draw [color=zzttqq] (-0.02,3.57)-- (0,-1);
\draw [color=zzttqq] (2,3.57)-- (2,-1);
\draw [color=zzttqq] (1,5)-- (1,3.57);
\end{tikzpicture}

%% file: circuit-figure.tex
\ctikzset{bipoles/length=.6cm}
\begin{circuitikz}[scale=0.6]
\draw[color=brown]
(0,0) to[short, *-*] (1,0)
      to[short, *-*] (2,0)
      to[short, *-] (4,0)
      to[sV=$\omega$](4,3)
      to[short, -*] (2,3)
      to[short, *-*] (0,3)
      to[C, *-*]     (0,2)
      to[L, *-*]     (0,0)
(1,0) to[L, *-*]     (1,2)
      to[short, *-*] (0,2)
(1,2) to[short, *-*] (2,2)        
(2,0) to[L, *-*]     (2,2)
      to[C, *-*]     (2,3)
;\end{circuitikz}